\documentclass[preprint,12pt]{elsarticle}
\usepackage{amsmath}
\usepackage{amsthm}
\usepackage{amssymb,epsfig}
\usepackage[cp 1250]{inputenc}
\usepackage{a4}
\usepackage{graphics}
\usepackage[dvips]{color}

\newtheorem{lemma}{Lemma}[section]
\newtheorem{prop}[lemma]{Proposition}
\newtheorem{coro}[lemma]{Corollary}
\newtheorem{thm}[lemma]{Theorem}
\theoremstyle{definition}
\newtheorem{definition}[lemma]{Definition}
\newtheorem{remark}[lemma]{Remark}
\newtheorem{example}[lemma]{Example}

\def\Z{\mathbb Z}

\def\A{\mathcal A}

\def\pf{\begin{proof}}
\def\pfk{\end{proof}}
\def\le{\leqslant}
\def\leq{\leqslant}
\def\ge{\geqslant}
\def\geq{\geqslant}

\begin{document}
\begin{frontmatter}
\title{Parallel addition\\ in non-standard numeration systems }
 
\author[liafa]{Christiane \textsc{Frougny}}
\author[doppler]{Edita \textsc{Pelantov\'a}\corref{cor1}}
\author[doppler]{Milena \textsc{Svobodov\'a}}

\cortext[cor1]{Corresponding author}
\address[liafa]{LIAFA, UMR 7089 CNRS \& Universit\'e Paris 7,  and Universit\'e Paris 8\\
Case 7014, 75205 Paris Cedex 13, France}
\address[doppler]{Doppler Institute for Mathematical Physics and
Applied Mathematics,\\ and
Department of Mathematics, FNSPE, Czech Technical University, \\
Trojanova 13, 120 00 Praha 2, Czech Republic}

\begin{abstract}
We consider numeration systems where digits are integers and the base is an algebraic number $\beta$ such that $|\beta|>1$ and $\beta$ satisfies a polynomial where one coefficient is dominant in a certain sense.
For this class of bases $\beta$, we can find an alphabet of signed-digits on which addition is realizable by a parallel
algorithm in constant time.
This algorithm is a kind of generalization of the one of Avizienis.
We also discuss the question of cardinality of the used alphabet, and we are able to modify our algorithm in order to work with a smaller alphabet.
We then prove that $\beta$ satisfies this dominance condition if and only if it has no
conjugate of modulus $1$.
When the base $\beta$ is the Golden Mean, we further refine the construction to obtain a parallel algorithm on the alphabet $\{-1,0,1\}$. This alphabet cannot be reduced any more.

\end{abstract}


\end{frontmatter}

\section{Introduction}

A positional numeration system is given by a base and by a set of digits.
The base $\beta$ is a real or complex number such that $|\beta|>1$, and the digit set $\mathcal{A}$
 is a finite alphabet of real or complex digits.
The most studied numeration systems are of course the usual ones, where the base is a positive integer.
But there have been also numerous studies where the base is an irrational real number (the so-called \emph{$\beta$-expansions}),
a complex number, or a non-integer rational number, \textit{etc}. Some surveys can be found in \cite[Chapter 7]{Lothaire} and \cite[Chapter 2]{cant}.

In this work, the base $\beta$ is an algebraic number, that is to say, the root of a polynomial with integer coefficients,
and the digits of $\mathcal{A}$ are consecutive integers. We consider only finite $\beta$-representations, see Section~\ref{prem} for definitions.
Let us denote
$$
{\rm Fin}_{\mathcal{A}}(\beta) = \Bigl\{ \ \sum_{i\in I}\
{x_i\beta^i} \ |  \ I \subset \mathbb{Z},  \ \ I \ \
\hbox{finite}, \ \ x_i \in \mathcal{A}\Bigr\}\,.
$$
In general, the set ${\rm Fin}_{\mathcal{A}}(\beta)$ is not closed under addition.
In this paper, we will prove that  when $\beta$ is a root of a polynomial with a dominant coefficient,
the alphabet $\mathcal{A}$ can be chosen in such a way that ${\rm Fin}_{\mathcal{A}}(\beta)$ is a ring, and,
 moreover, addition can be performed by a parallel algorithm with time complexity $\mathcal{O}(1)$.

Addition of two numbers in the classical $b$-ary numeration
system, where $b$ is an integer $\ge 2$, has linear time
complexity. In order to save time,
 several solutions have been proposed.
A popular one is the Avizienis signed-digit
representation~\cite{Avizienis},
 which consists
of changing the digit set. Instead of taking digits from the
canonical alphabet $\{0, \ldots, b-1\}$, they are taken from a
symmetric alphabet of the form $\{-a,\ldots,0, \ldots,a\}$, $a$
being an integer such that
 $b/2  < a \le b-1$ ($b$ has to be $\ge 3$).
 Such a numeration system is {\em redundant}, that is to say, some
 numbers may have
several representations.  This property allows one to perform
addition in constant time in parallel, because there is a limited
carry propagation. A similar algorithm for base 2 has been devised
by Chow and Robertson~\cite{ChowRobertson}, using alphabet $\{-1,
0,1\}$. Here addition is realized in parallel with a window of
size 3. Notice that Cauchy~\cite{Cauchy} already considered the
redundant representation in base 10 and alphabet $\{-5,\ldots,0,
\ldots,5\}$.

In symbolic dynamics, such functions computable in parallel
are called {\em local}, more precisely $p$-local, which means that
to determine the image of a word by a $p$-local function, it is
enough to consider a sliding window of length $p$ of the input. It
is a general result that a $p$-local function is computable by an
on-line finite automaton with delay $p-1$, see~\cite{Frougny2} for
instance.

\medskip

Amongst the non-standard bases, special attention has been paid to
the complex ones, which allow to represent any complex number by a
single sequence (finite or infinite) of natural digits, without
separating the real and the imaginary part. For instance, it is
known that every complex number can be expressed with base $-1+i$
and digit set $\{0,1\}$, \cite{Penney}.

Parallel algorithms for addition in bases $-2$, $i\sqrt{2}$, $2i$
and $-1+i$ have been given in \cite{NM}. Results on addition in
bases $-b$, $i \sqrt{b}$ and $-1+i$ in connection with automata
theory have been presented in \cite{Frougny2}. In particular, if
$b$ is an integer, $|b| \ge 3$, and $\A=\{-a, \ldots,0,\ldots a\}$
with $a=\lceil \tfrac{|b|+1}{2} \rceil$, addition in base
$\beta=\sqrt[q]{b}$ is computable in parallel and is a
$(q+1)$-local function. If $|b| \ge 2$ is even, set $a=|b| /2$;
then addition in base $\beta=\sqrt[q]{b}$ on $\A$ is a
$(2q+1)$-local function.

\bigskip

The paper is organized as follows: First we suppose that the base $\beta$ is a root of a polynomial with integer coefficients such that one of them is greater than twice the sum of the moduli of the other ones.
We then say that $\beta$ satisfies the {\em strong representation of zero property} (or, for short, that {\em $\beta$ is SRZ}).
We give a parallel algorithm, Algorithm~I, for addition in that case, Theorem~\ref{AI}.
When $\beta$ is an integer $\ge 3$, it is the same algorithm as the one of Avizienis.

Section~\ref{reduc} is devoted to reduction of the working alphabet. We now suppose that $\beta$ is a root of a polynomial such that one coefficient is greater than the sum of the moduli of the other ones.
We then say that $\beta$ satisfies the {\em weak representation of zero property} (or, for short, that {\em $\beta$ is WRZ}).
We give a parallel algorithm, Algorithm~II, where the alphabet is reduced compared to Algorithm~I, but there is a fixed number of iterations depending (only) on the weak polynomial satisfied by $\beta$, to be  compared with Algorithm~I which always has just one iteration.

We then show that in fact all algebraic numbers with no conjugate of modulus~$1$ are SRZ (or WRZ), and we give a constructive method to obtain the suitable polynomial (strong or weak representation of zero) from the minimal polynomial of $\beta$, Proposition \ref{dominant}.

In the end, we concentrate on the case where $\beta$ is the Golden Mean.
Algorithm~II works with alphabet  $\{-3, \ldots,3\}$, which is big compared to the minimally redundant alphabet.
Using ideas similar to the Chow and Robertson algorithm, we obtain a parallel algorithm on $\{-1,0,1\}$, Algorithm~III.
This algorithm can be applied to the Fibonacci numeration system as well.

%%%%%%%%%%%%%%%%%%%%%%%%%%%%%%%%%%%%%%%%%%%%%%%%%%%%%%%%%%%%%%%%%%%%%%%%%%%%%%%%%%%%%%%%%
%%%%%%%%%%%%%%%%%%%%%%%%%%%%%%%%%%%%%%%%%%%%%%%%%%%%%%%%%%%%%%%%%%%%%%%%%%%%%%%%%%%%%%%%%

\section{Preliminaries}\label{prem}

A finite {\em word} on the alphabet $\mathcal{A}$ is a concatenation of a finite number of letters of $\mathcal{A}$. The set of words on $\mathcal{A}$ is the free monoid $\mathcal{A}^*$. The set of bi-infinite words on $\mathcal{A}$ is denoted $\A^\Z$.

A finite $\beta$-{\em representation} of $x$ with digits in $\mathcal{A}$ is a finite sequence $(x_i)_{m \le i \le n}$,  with $x_i$ in $\mathcal{A}$, such that $x=\sum_{i=m}^{n} x_i \beta^i$. It is denoted by the word $$x_{n}  x_{n-1} \cdots x_{m+1} x_{m}$$ with the most significant digit at the left hand-side, as in the decimal numeration system.

The notion of a function computable in parallel comes from computer arithmetic (see \cite{Kornerup}), where it is defined on finite words, but we give here a definition on bi-infinite words.
A function $\varphi : \mathcal{A}^{\Z} \rightarrow \mathcal{B}^{\Z}$ is said to be {\em computable in parallel}
if there exist two non-negative integers $r$ and $t$, a positive integer $k$, and a function $\Phi$
from $\mathcal{A}^p$ to $\mathcal{B}^k$, with $p=r+t+k$, such that if $u=(u_i)_{i \in \Z} \in \mathcal{A}^{\Z}$ and $v=(v_i)_{i\in \Z}
\in \mathcal{B}^{\Z}$, then $v=\varphi(u)$ if, and only if, for every $i$ in $\Z$, $v_{ki+k-1}\cdots v_{ki}=\Phi(u_{ki+k+t-1} \cdots u_{ki-r})$\footnote{Careful! Indices of $\Z$ are decreasing from left to right.}.
This means that the image of $u$ by $\varphi$ is obtained through a sliding window of length $p$.
Such functions are computable in constant time in parallel.

The notion of a {\em local function} comes from symbolic dynamics (see \cite{LM}) and is often called a {\em sliding block code}.
It is a function computable in parallel with $k=1$.
The parameter $r$ is called the \emph{memory} and the parameter $t$ is called the \emph{anticipation} of the function $\varphi$.
The function $\varphi$ is then said to be {\em $p$-local}.

\bigskip

To be self-contained, we recall the classical algorithms for parallel addition of Avizienis~\cite{Avizienis}, and of Chow and
Robertson~\cite{ChowRobertson}. In what follows, ``for each $i$ in parallel"
means that ``each numbered step is executed in parallel, and the results of the parallel
computations are shared between the steps".

\vskip1.5cm

\hrule

\vskip0.2cm

\noindent {\bf Algorithm of Avizienis}: Base $\beta=b$, $b\ge 3$
integer, parallel addition on alphabet $\A=\{-a,\ldots,0, \ldots,a\}$, $b/2  < a \le b-1$.

\vskip0.2cm

\hrule

\vskip0.2cm

\noindent{\sl Input}:  two words $x_n \cdots x_m$ and $y_n \cdots
y_m$ of $\mathcal{A}^*$, with $m \le n$, $ x = \sum_{i=m}^{n}
x_i\beta^i$ and
$y = \sum_{i=m}^{n} y_i\beta^i$.\\
{\sl Output}: a word $z_{n+1} \cdots z_{m}$ of $\mathcal{A}^*$
such that
$$z =x+y = \sum_{i=m}^{n+1} z_i\beta^i.$$

\noindent\texttt{for each $i$ in parallel do}\\
0. \hspace*{0.5cm} $z_i := x_i+y_i$\\
1. \hspace*{0.5cm} \texttt{if}  $z_i  \ge a$ \texttt{then} $q_{i}:=1$, $r_i:=z_i-b$\\
\hspace*{1cm} \texttt{if}  $z_i  \le -a$ \texttt{then} $q_{i}:=-1$, $r_i:=z_i+b$\\
\hspace*{1cm} \texttt{if}  $-a+1 \le z_i  \le a-1$ \texttt{then} $q_{i}:=0$, $r_i:=z_i$\\
2. \hspace*{0.5cm} $z_i:= q_{i-1} + r_i$

\vskip0.2cm

\hrule \vskip0.2cm Addition realized by the Avizienis algorithm is a $2$-local function,
with memory $1$ and anticipation $0$.
Notice that the minimally redundant symmetric alphabet is obtained with the value $a=\lceil \tfrac{b+1}{2} \rceil$.

\bigskip

The Chow and Robertson algorithm works for base $2$ and alphabet $\{-1,0,1\}$. We give here a generalization to an even base
$\beta=b=2a$, $\A=\{-a,\ldots,0, \ldots,a\}$.

\vskip0.2cm

\hrule

\vskip0.2cm

\noindent {\bf Algorithm of Chow and Robertson}: Base $\beta = b =2a$, $a\ge 1$ integer, parallel addition on alphabet $\A=\{-a,\ldots,0, \ldots,a\}$.

\vskip0.2cm

\hrule

\vskip0.2cm

\noindent{\sl Input}:  two words $x_n \cdots x_m$ and $y_n \cdots
y_m$ of $\mathcal{A}^*$, with $m \le n$, $ x = \sum_{i=m}^{n}
x_i\beta^i$ and
$y = \sum_{i=m}^{n} y_i\beta^i$.\\
{\sl Output}: a word $z_{n+1} \cdots z_{m}$ of $\mathcal{A}^*$
such that
$$z =x+y = \sum_{i=m}^{n+1} z_i\beta^i.$$

\noindent\texttt{for each $i$ in parallel do}\\
0. \hspace*{0.5cm} $z_i := x_i+y_i$\\
1. \hspace*{0.5cm} \texttt{if}  $a+1 \le z_i  \le b$ \texttt{then} $q_{i}:=1$, $r_i:=z_i-b$\\
\hspace*{1cm} \texttt{if}  $-b \le z_i  \le -a-1$ \texttt{then} $q_{i}:=-1$, $r_i:=z_i+b$\\
\hspace*{1cm} \texttt{if}  $-a+1 \le z_i  \le a-1$ \texttt{then} $q_{i}:=0$, $r_i:=z_i$\\
\hspace*{1cm} \texttt{if}  $z_i  = a$ \texttt{and} $z_{i-1} >0$ \texttt{then} $q_{i}:=1$, $r_i:=-a$\\
\hspace*{1cm} \texttt{if}  $z_i  = a$ \texttt{and} $z_{i-1} \le 0$ \texttt{then} $q_{i}:=0$, $r_i:=a$\\
\hspace*{1cm} \texttt{if}  $z_i  = -a$ \texttt{and} $z_{i-1} <0$ \texttt{then} $q_{i}:=-1$, $r_i:=a$\\
\hspace*{1cm} \texttt{if}  $z_i  = -a$ \texttt{and} $z_{i-1} \ge 0$ \texttt{then} $q_{i}:=0$, $r_i:=-a$\\
2. \hspace*{0.5cm} $z_i:= q_{i-1}  + r_i$

\vskip0.2cm

\hrule \vskip0.2cm Addition realized by the Chow and Robertson algorithm is a $3$-local function, with memory $2$ and anticipation $0$.

The main difference between the two algorithms is that the Avizienis algorithm is \emph{neighbour-free}, while the Chow and Robertson algorithm is \emph{neighbour-sensitive}, since the decision taken at position $i$ in Step 1 depends also on the digit at position $i-1$.

%%%%%%%%%%%%%%%%%%%%%%%%%%%%%%%%%%%%%%%%%%%%%%%%%%%%%%%%%%%%%%%%%%%%%%%%%%%%%%%%%%%%%%%%%
%%%%%%%%%%%%%%%%%%%%%%%%%%%%%%%%%%%%%%%%%%%%%%%%%%%%%%%%%%%%%%%%%%%%%%%%%%%%%%%%%%%%%%%%%

\section{Parallel algorithm}

Now we move to bases $\beta$ that are algebraic numbers. Let us formalize the assumption we want $\beta$ to satisfy.

\begin{definition}\label{DefSRZ}
Let $\beta$ be such that $|\beta|>1$. We say that $\beta$ satisfies the {\em strong representation of zero property}
(or, for short, that {\em $\beta$ is SRZ}) if there exist integers  $b_k,  b_{k-1}, \ldots,  b_1, b_0, b_{-1},
\ldots,b_{-h}$, for some non-negative integers $h$ and $k$, such that $\beta$ is a root of the polynomial
\begin{equation}\label{treba1}
S(X)=b_kX^k + b_{k-1}X^{k-1} + \cdots  + b_1X + b_0 + b_{-1}X^{-1} + \cdots + b_{-h}X^{-h}
\end{equation}
and
$$ b_0 > 2 \sum_{i\in \{-h, \ldots,k\}\setminus \{0\}}{|b_i|}\,.$$
The polynomial $S$ is a said to be a {\em strong polynomial} for $\beta$.
\end{definition}

If $\beta$ satisfies (\ref{treba1}), then the word $b_k b_{k-1}\cdots b_1 b_0 b_{-1} \cdots b_{-h}$ is a $\beta$-representation of zero.
From this we can derive a set of rewriting rules, generated by the rule $b_k b_{k-1}\cdots b_1 b_0 b_{-1} \cdots b_{-h} \mapsto 0$.

To simplify the notation, we set $ B= b_0$ and $ M = \sum\limits_{i\in \{-h, \ldots,k\}\setminus \{0\}}{|b_i|}\,.$ The inequality from Definition \ref{DefSRZ} now reads
\begin{equation}\label{treba2}
B > 2 M\,.
\end{equation}
Suppose that $\beta$ is SRZ. We choose the symmetric alphabet
\begin{equation}\label{treba3}
\mathcal{A} = \{-a, \ldots,  0, \ldots,  a\}, \  \hbox{ where} \
a=\left\lceil\tfrac{B-1}{2}\right\rceil + \left\lceil\tfrac{B-1}{2(B-2M)}\right\rceil M\,.
\end{equation}
For this fixed alphabet $\mathcal{A}$, we describe a parallel algorithm for addition in base $\beta$.
Let us denote \begin{equation}\label{treba4} c = \left\lceil\tfrac{B-1}{2(B-2M)}\right\rceil \quad \hbox{and} \quad
a' = \left\lceil\tfrac{B-1}{2}\right\rceil\,.\end{equation} Then $a=a'+cM$. The alphabet $\mathcal{A}' = \{ -a', \ldots, 0,
\ldots, a'\} \subset \mathcal{A}$ will be referred to as the \emph{inner} alphabet.

\bigskip

\hrule

\vskip0.2cm

\noindent {\bf Algorithm~I}: Parallel addition for base $\beta$ with the strong representation of zero property ($\beta$ is SRZ).

\vskip0.2cm

\hrule

\vskip0.2cm

\noindent{\sl Input}:  two words $x_n \cdots x_m$ and $y_n \cdots
y_m$ of $\mathcal{A}^*$, with $m \le n$, $ x = \sum_{i=m}^{n}
x_i\beta^i$ and
$y = \sum_{i=m}^{n} y_i\beta^i$.\\
{\sl Output}: a word $z_{n+k} \cdots z_{m-h}$ of $\mathcal{A}^*$
such that
$$z =x+y = \sum_{i=m-h}^{n+k} z_i\beta^i.$$

\noindent\texttt{for each $i$ in parallel do}\\
0. \hspace*{0.5cm} $z_i := x_i+y_i$\\
1. \hspace*{0.5cm} find $q_i \in \{-c, \ldots, 0, \ldots, c\}$
such that $z_i -q_iB \in \mathcal{A}'$\\
2. \hspace*{0.5cm} $z_i:= z_i  - \sum\limits_{j=-h}^kq_{i-j}b_j$

\vskip0.2cm

\hrule

\vskip0.5cm

Before proving the correctness of the algorithm, let us stress that the numbers $c$, $a$, and $a'$ have been defined in such a way that
 they satisfy
\begin{equation}\label{nerovnosti}
(i)\quad 2a'+1 \geq B, \qquad   (ii)\quad a'+cM \leq a,  \qquad (iii)\quad 2a-cB\leq a'.
\end{equation}
The first two inequalities   are clear.  To verify the last one, we will use the fact that, for any positive real numbers $\gamma, \delta$,
one has $\lceil \gamma \delta\rceil \leq \lceil \gamma\rceil \lceil \delta\rceil$. Putting $\gamma=\frac{B-1}{2(B-2M)}$ and
$\delta=B-2M$, we obtain $$ a' = \left\lceil\tfrac{B-1}{2}\right\rceil = \left\lceil\tfrac{B-1}{2(B-2M)}
(B-2M)\right\rceil\leq \left\lceil\tfrac{B-1}{2(B-2M)}
\right\rceil (B-2M) = c(B-2M)\,.$$ This inequality, together with our choice $a = a'+cM$, already proves the third inequality of
\eqref{nerovnosti}.

\begin{thm}\label{AI}
Suppose that $\beta$ is SRZ. Then Algorithm~I realizes addition in constant time in parallel in ${\rm Fin}_{\mathcal{A}}(\beta)$
 with $\mathcal{A}= \{-a, \ldots,  0, \ldots,  a\}$, where $a=\left\lceil\tfrac{B-1}{2}\right\rceil + \left\lceil\tfrac{B-1}{2(B-2M)}\right\rceil M$.
\end{thm}

\pf The proof is structured in two parts:

\begin{itemize}

    \item
\noindent The digits of the output belong to $\mathcal{A}$:

\noindent After Line 0 of Algorithm~I,  $z_i$ belongs to $\{- 2a,
\ldots, 0, \dots, 2a\}$.\\
Since the inner alphabet $\mathcal{A}'$ is formed by $2a'+1$ consecutive integers and $B\leq 2a'+1$, there exists
$q_i$ such that  $z_i -q_iB \in \mathcal{A}'$. Since $|z_i | \leq 2a$ and $2a-cB\leq a'$,
the integer $q_i$ can be found in range $  \{-c, \ldots, 0, \ldots, c\}$.\\
Finally, in Line 2, the new digit $z_i$  (denoted for the moment by $z_i^{new}$) is
 $$z_i^{new}= \underbrace{z_i  -  q_i B}_{:= E} -
 \underbrace{\sum\limits_{j=-h}^{-1}q_{i-j}b_j-\sum\limits_{j=1}^{k}q_{i-j}b_j}_{:= F}.$$
As $|E| \leq a'$ and $|F| \leq cM$, the resulting $z_i^{new}$ satisfies $|z_i^{new}|\leq a'+cM \leq  a$.

    \item
\noindent The output represents the number $x+y$:

\noindent In order to avoid a tedious description  of ranges of summation indices, we will consider all coefficients which do not play any
role in our  consideration to be equal to zero.
With this convention, $b_i=0$ for all $i\notin \{-h, \ldots, 0,\ldots, k\}$, and therefore we may express the rewriting rule \eqref{treba1} as
 $\sum\limits_{i\in \mathbb{Z}}b_i\beta^i = 0$.
Similarly,  $x_i$ and $y_i$ are set to $0$ for all $i\notin \{m, \ldots, n\}$, and thus $x+y = \sum\limits_{i\in \mathbb{Z}}(x_i+y_i)\beta^i$.
Also the auxiliary coefficients $q_i$ are set to $0$ for all $i \notin \{m, \ldots, n\}$.\\

After Line 2 of Algorithm~I, we obtain
 $$z=\sum_{i\in
   \mathbb{Z}}z_i^{new}\beta^i =
   %\sum_{i\in
  % \mathbb{Z}}w_i\beta^i -  \sum_{i\in
   %\mathbb{Z}} \Bigl(\sum_{j\in
   %\mathbb{Z}}q_{i-j}b_j\Bigr)\beta^i=
   x+y - \sum_{i\in
   \mathbb{Z}} \Bigl(\sum_{j\in
   \mathbb{Z}}q_{i-j}b_j\Bigr)\beta^i\,.$$
Since for the last sum  we have
$$\sum_{i\in
   \mathbb{Z}} \Bigl(\sum_{j\in
   \mathbb{Z}}q_{i-j}b_j\Bigr)\beta^i = \sum_{j\in
   \mathbb{Z}}\sum_{\ell\in
   \mathbb{Z}} q_\ell b_j\beta^{\ell+j} =  \sum_{\ell\in
   \mathbb{Z}}q_\ell \beta^\ell\underbrace{\sum_{j\in
   \mathbb{Z}} b_j\beta^{j}}_{=0} = 0 \,,$$
the output is correct.
\end{itemize}\pfk

\begin{remark}
The correctness of the algorithm stems from the inequalities \eqref{nerovnosti}, which are satisfied by many triplets  $a$, $a'$, $c$.
Our specific choice of $a$, $a'$, and $c$ displayed in \eqref{treba3} and \eqref{treba4} gives the smallest possible value $a$,
 and thus the smallest cardinality of the alphabet $\mathcal{A}$ in our algorithm.
Of course, this does not exclude existence of other parallel algorithms with smaller alphabets.
Let us demonstrate the optimality of our choice  of $a$, $a'$, and $c$.
According to (ii) of \eqref{nerovnosti}, in order to minimize $a$, we have to minimize integers $a'$ and $c$.
The choice of $a'$ in \eqref{treba4} is minimal with respect to the requirement $2a'+1\geq B$.
By (ii) and (iii) of \eqref{nerovnosti}, we obtain
$$ 2a'+1 \leq 2a -2cM +1 \leq a'+cB-2cM+1 \quad \Longrightarrow \quad a'\leq cB-2cM\,. $$
Combining the last inequality with (i) of \eqref{nerovnosti}, and insisting on positivity of $c$, we have
$$ B\leq 2a'+1 \leq 2cB - 4cM + 1 \quad \Longrightarrow \quad c\geq \tfrac{B-1}{2(B-2M)},$$
 which already shows that our choice of $c$ in \eqref{treba4} was minimal as well.
Let us point out that we  also used  the inequality $B>2M$.
\end{remark}

\begin{coro}
If $\beta$ is SRZ, then addition realized by Algorithm~I is a $(h+k+1)$-local function with memory $k$ and anticipation $h$. Algorithm~I is neighbour-free.
\end{coro}

\begin{example} Consider  the base $\beta = 10$. It is SRZ for the polynomial $-X+10=0$, where $B =10$ and $M=1$. In this case
$$c=\left\lceil\tfrac{9}{16}\right\rceil  = 1, \quad  a' =
\left\lceil\tfrac{9}{2}\right\rceil = 5, \quad \hbox{and}\quad a=
6\,.$$
Therefore, our algorithm provides parallel addition in the decimal numeration system on alphabet $\mathcal{A} = \{-6, \ldots, 0, \ldots, 6\}$,
 and in fact it is precisely the algorithm that Avizienis gave in 1961, see \cite{Avizienis}.

More generally, for any integer $b\geq 3$, the base $\beta = b$ is SRZ for the polynomial $-X+b=0$, and therefore
 we can apply the same algorithm for addition, namely on the alphabet $\{-a, \ldots, 0,\ldots , a\}$ with $a=\left\lceil\tfrac{b+1}{2}\right\rceil$.
\end{example}

\begin{example}\label{binar} Consider  the base $\beta = 2$. For such a base, $-X+2=0$ is not a strong polynomial.
Nevertheless, this base satisfies also the polynomial $-X^2 + 4=0$, which already is strong, with  $B =4$ and $M=1$.
Now we have
$$c=\left\lceil\tfrac{3}{4}\right\rceil  = 1, \quad  a' = \left\lceil\tfrac{3}{2}\right\rceil = 2, \quad \hbox{and}\quad a= 3\,.$$
So Algorithm~I works for base $2$ with the alphabet $\{ -3,\ldots,0,\ldots,3\}$, and is $3$-local.
Remind that the Chow and Robertson algorithm is $3$-local as well, but it works with smaller alphabet $\{-1, 0, 1\}$.
\end{example}

\begin{example}\label{fibonacci} Let us consider the base  $\beta = \frac{1+\sqrt{5}}{2}$, the Golden Mean.
This base $\beta$ is one root of the quadratic equation $X^2=X+1$, the second root is $\beta' = \frac{1-\sqrt{5}}{2} = -\frac{1}{\beta}$.
Since $\beta^4+(\beta')^4 = 7$, $\beta$ is a root of the strong polynomial $$ S(X)=-X^4 + 7 - \tfrac{1}{X^4}$$
with $B=7$ and $M=2$. This implies $c=1$, $a'=3$, and $a=5$.
The alphabet we work with is $\mathcal{A} = \{ -5, \ldots, 0,\ldots, 5\}$.
The inner alphabet is $\mathcal{A'} = \{ -3, \ldots,0,  \ldots, 3\}$.

In the following, we denote the signed-digit
$-b$ by $\overline{b}$ (this notation is already present in the work of Cauchy~\cite{Cauchy}).
Below is shown the performance of Algorithm I for addition of two numbers from ${\rm Fin}_{\mathcal{A}}(\beta)$, namely
$$ x = 2\beta^7 +5\beta^6 -2\beta^5 +5\beta^4 - 5\beta^3 +3\quad  \hbox{and} \quad  y = 5\beta^8 +\beta^7 + 2\beta^6 -2\beta^5 +5\beta^4-4\beta^3 +5.$$
We have $x_7 \cdots x_0=25\overline{2}5\overline{5}003$ and $y_8 \cdots y_0=512\overline{2}5\overline{4}005$.

$$\begin{array}{llllllllllcllllllll}
x&\mapsto~~~~~ & &&&&&2&5&\overline{2}&5&\overline{5}&0&0&3&&&&\\
y&\mapsto~~~~~ & &&&&5&1&2&\overline{2}&5&\overline{4}&0&0&5&&&&\\
 \hline
z&\mapsto~~~~~ & &&&&5&3&7&\overline{4}&10&\overline{9}&0&0&8&&&&\\
\end{array}$$
All positions $i$ of $z$  for which $z_i$ does not belong to the inner alphabet  $\{-3, \ldots, 0,  \ldots,3\}$ must be modified by adding or
subtracting the rewriting rule $\overline{1} 0007000\overline{1}$, with the digit $7$ on position $i$, as follows:

$$\begin{array}{llllllllllcllllllll}
z&\mapsto~~~~~ & &&&&5&3&7&\overline{4}&10&\overline{9}&0&0&8&&&&\\
0&\mapsto~~~~~ &1 &0&0&0&\overline{7}&0&0&0&1&&&&&&&\\
0&\mapsto~~~~~ &&&1 &0&0&0&\overline{7}&0&0&0&1&&&&&\\
0&\mapsto~~~~~ &&&&\overline{1} &0&0&0&7&0&0&0&\overline{1}&&&&\\
0&\mapsto~~~~~ &&&&&1 &0&0&0&\overline{7}&0&0&0&1&&&\\
0&\mapsto~~~~~ &&&&&&\overline{1} &0&0&0&7&0&0&0&\overline{1}&&\\
0&\mapsto~~~~~ &&&&&&&&&1 &0&0&0&\overline{7}&0&0&0&1\\

\hline

z&\mapsto~~~~~ &1 &0&1&\overline{1}&\overline{1}&2&0&3&5&\overline{2}&1&\overline{1}&2&\overline{1}&0&0&1\\
\end{array}
$$
Hence $x+y=z=\beta^{12} +
\beta^{10}-\beta^9-\beta^8+2\beta^7+3\beta^5+5\beta^4-2\beta^3
+\beta^2 -\beta +2-\beta^{-1}+\beta^{-4}$.
\end{example}

\begin{example}\label{Gaussian} Consider the complex base $\beta = -1+i$.
Since $\beta^4 = -4$, we can use the strong polynomial $X^4 +4$ with $B=4$ and $M=1$, \textit{i.e.},
 we can perform parallel addition  on the  alphabet $\mathcal{A}=\{-3, \ldots, 0,  \ldots,3\}$, and this addition is a $5$-local function.
This result appeared in \cite{Frougny2}.
There is also a parallel algorithm on alphabet $\{-2,\ldots,2\}$ (\`a la Chow and Robertson), which gives a $9$-local function, \cite{Frougny2}.
On the minimally redundant alphabet $\{-1,0,1\}$, there is a parallel addition algorithm provided by Nielsen and Muller~\cite{NM}.
\end{example}

\begin{example}\label{rat}
Take for base $\beta$ the rational $\tfrac{7}{2}$.
This is an algebraic number, which is not an algebraic integer.
It is known that any non-negative integer has a finite expansion in this base on the alphabet $\{0,\ldots,6\}$, see~\cite{FrougnyKlouda}.
There is a strong polynomial $S(X)=-2X+7$ with $B=7$ and $M=2$. Thus $a'=3$, $c=1$, and $a=5$.
It is not difficult to see that any integer has a finite representation on the redundant alphabet $\{-5, \ldots,0,\ldots,5\}$,
 and, by Algorithm~I, addition is realizable in parallel.
\end{example}

%%%%%%%%%%%%%%%%%%%%%%%%%%%%%%%%%%%%%%%%%%%%%%%%%%%%%%%%%%%%%%%%%%%%%%%%%%%%%%%%%%%%%%%%%
%%%%%%%%%%%%%%%%%%%%%%%%%%%%%%%%%%%%%%%%%%%%%%%%%%%%%%%%%%%%%%%%%%%%%%%%%%%%%%%%%%%%%%%%%

\section{Reduction of the  alphabet}\label{reduc}

It is very difficult to prove that an alphabet used by some algorithm for parallel addition is minimal.
It is a lot easier to prove that the alphabet is not minimal; which is the case for the alphabet we worked with in Algorithm I.
Let us start this section with  an example.

\begin{example}\label{golden} Take for $\beta$ the Golden Mean. This base satisfies the equation
$$-\beta^2 + 3 - \frac{1}{\beta^2} = 0 \,.$$
We say that rewriting rule $\bar1 0 3 0\bar1 \mapsto 0$ is \emph{positively applied} at position $i$
if the number 3 is added to the actual digit occurring at position $i$ and the number $1$ is subtracted from the actual digits at
positions $i-2$ and $i+2$.
Analogously, we define the \emph{negatively applied} at position $i$.
We shall present a parallel algorithm for addition on the alphabet
$$\mathcal{A} = \{-3,\ldots, 0, \ldots, 3\}$$

Let us explain the algorithm less formally.

%\vskip0.2cm

\bigskip

\hrule

\vskip0.2cm

\noindent{\sl Input}:  two finite sequences of digits $(x_i)$ and
$(y_i)$ of $\{-3,\ldots,3\}$, with $ x = \sum x_i\beta^i$ and
$y = \sum y_i\beta^i$.\\
{\sl Output}: a finite sequence of digits $(z_i)$ of
$\{-3,\ldots,3\}$ such that
$$z =x+y = \sum z_i\beta^i.$$

\noindent\texttt{for each $i$ in parallel do}\\
$\begin{array}{lll}
\hspace*{-0.15cm}0. \hspace*{0.5cm}  & z_i := x_i+y_i &   \\
\hspace*{-0.15cm}1.  \hspace*{0.5cm} & \texttt{if}  \; z_i  \in \{2,3,4,5,6\} & \texttt{then} \;
    \hbox{ apply negatively the rule at $i$} \\
& \texttt{if} \;  z_i  \in \{-6, -5, -4, -3, -2\}  & \texttt{then} \;
    \hbox{ apply positively the rule at $i$ }\\
& \texttt{if}  \;  z_i  \in \{-1,  0, 1\} & \texttt{then} \; \hbox{  do  nothing }\\
\hspace*{-0.15cm}2.  \hspace*{0.5cm}  &\texttt{if}  \; z_i  \in \{2,3,4,5\} & \texttt{then}  \;
    \hbox{ apply negatively the rule at $i$}\\
 &\texttt{if}  \;  z_i  \in \{-5, -4, -3, -2\} & \texttt{then}  \;
    \hbox{ apply positively the rule at $i$}\\
& \texttt{if}   \;  z_i  \in \{-1, 0,1 \} & \texttt{then} \; \hbox{  do  nothing }\\
\hspace*{-0.15cm}3.  \hspace*{0.5cm}  &\texttt{if}  \;    z_i  \in \{2,3,4\} & \texttt{then} \;
    \hbox{ apply negatively the rule at $i$}\\
& \texttt{if}  \;  z_i \in \{-4, -3, -2\} & \texttt{then}  \;
    \hbox{ apply positively the rule at $i$}\\
& \texttt{if}  \;  z_i  \in \{-1, 0,1 \} & \texttt{then}  \; \hbox{  do  nothing }\\
  \end{array}$\\

\vskip0.2cm

\hrule

\vskip0.5cm

\noindent During the course of the algorithm, we change only the actual values of the coefficients $z_i$, but we do not change the sum  $\sum z_i\beta^i$.
Therefore, for correctness of the algorithm, we just have to realize that
\begin{itemize}
    \item after Step 0, each $z_i$ belongs to $ \{ -6, \ldots, 6\}$,
    \item after Step 1, each $z_i$ belongs to $ \{ -5, \ldots, 5\}$,
    \item after Step 2, each $z_i$ belongs to $ \{ -4, \ldots, 4\}$, and, finally,
    \item after Step 3, each $z_i$ belongs already to the original alphabet $\mathcal{A} = \{ -3, \ldots, 3 \}$.
\end{itemize}

Let us illustrate the algorithm on the following example:
$$\begin{array}{lllllllllllcllllllll}
\hline
\mathrm{Step\ 0.\ \ } & x&\mapsto~~~~~ & &&&&3&\overline{1}&3&0&3&&&&&&&&\\
&y&\mapsto~~~~~ & &&&&2&0&3&\overline{2}&3&&&&&&&&\\
\hline
&z&\mapsto~~~~~ & &&&&5&\overline{1}&6&\overline{2}&6&&&&&&&&\\
\hline
\mathrm{Step\ 1.} & 0&\mapsto~~~~~ & &&1&0&\overline{3}&0&1&&&&&&&&&&\\
&0&\mapsto~~~~~ & &&&&1&0&\overline{3}&0&1&&&&&&&&\\
&0&\mapsto~~~~~ && &&&&\overline{1}&0&3&0&\overline{1}&&&&&&&\\
&0&\mapsto~~~~~ & && &&&&1&0&\overline{3}&0&1&&&&&&\\
\hline
&z&\mapsto~~~~~& &&1&0&3&\overline{2}&5&1&4&\overline{1}&1&&&&&&\\
\hline
\mathrm{Step\ 2.} &0&\mapsto~~~~~ & &&1&0&\overline{3}&0&1&&&&&&&&&&\\
&0&\mapsto~~~~~  &&&&\overline{1}&0&3&0&\overline{1}&&&&&&&&&\\
&0&\mapsto~~~~~ & &&&&1&0&\overline{3}&0&1&&&&&&&&\\
&0&\mapsto~~~~~ & && &&&&1&0&\overline{3}&0&1&&&&&&\\
\hline
&z&\mapsto~~~~~& &&2&\overline{1}&1&1&4&0&2&\overline{1}&2&&&&&&\\
\hline
\mathrm{Step\ 3.} &0&\mapsto~~~~~ &1&0&\overline{3}&0&1&&&&&&&&&&&&\\
&0&\mapsto~~~~~ & &&&&1&0&\overline{3}&0&1&&&&&&&&\\
&0&\mapsto~~~~~ & &&&&&&1&0&\overline{3}&0&1&&&&&&\\
&0&\mapsto~~~~~ & &&&&&&&&1&0&\overline{3}&0&1&&&&\\
\hline
&z&\mapsto~~~~~ &1&0&\overline{1}&\overline{1}&3&1&2&0&1&\overline{1}&0&0&1&&&&\\
\hline\hline
\end{array}
$$
\end{example}

When  $\beta$ is  the Golden Mean,  we have  seen   how using  a
weaker representation of zero $-\beta^2 + 3 - \frac{1}{\beta^2} = 0$
enables to exploit   the reduced alphabet $\{-3, \ldots, 3\}$
instead of the alphabet $\{-5, \ldots, 5\}$, which was necessary
for applying the strong representation of zero $-\beta^4 + 7  -
\frac{1}{\beta^4} = 0$. The idea of this reduction  can be
generalized to other bases $\beta$  as well.

\begin{definition}
Let $\beta$ be such that $|\beta|>1$. We say that $\beta$ satisfies the {\em weak representation of zero property}
(or, for short, that {\em $\beta$ is WRZ}) if there exist integers  $b_k,  b_{k-1}, \ldots,  b_1, b_0, b_{-1},
\ldots,b_{-h}$, for some non-negative integers $h$ and $k$, such that $\beta$ is a root of the polynomial
\begin{equation}
W(X)=b_kX^k + b_{k-1}X^{k-1} + \cdots  + b_1X + b_0 + b_{-1}X^{-1} + \cdots + b_{-h}X^{-h}
\end{equation}
and
\begin{equation}\label{trebaII1}
b_0 >  \sum_{i\in \{-h, \ldots,k\}\setminus \{0\}}{|b_i|}\,.
\end{equation}
The polynomial $W$ is said to be a \emph{weak polynomial} for $\beta$.
\end{definition}

When $\beta$ is WRZ, we can describe a parallel  algorithm for addition. Let us put
as above $M=\sum_{i\in \{-h, \ldots,k\}\setminus \{0\}}{|b_i|}$, and let
\begin{equation}\label{trebaII2}  \mathcal{A} =
\{ -a,  \ldots,  0,  \ldots, a\}, \ \  \hbox{where} \ \ a =
\left\lceil\tfrac{B-1}{2}\right\rceil +M\,.
\end{equation}
Similarly to Algorithm~I, the inner alphabet is $\mathcal{A}' = \{ -a', \ldots,  0,  \ldots, a'\}$ with $a'=\lceil\tfrac{B-1}{2}\rceil$.
The algorithm  works in
\begin{equation}\label{trebaII3} s+1\  \ \hbox{steps, where} \  \ s =
\left\lceil\tfrac{a}{B-M}\right\rceil \,.
\end{equation} The steps will be indexed by
 $\ell   = 0, 1, \ldots, s$.  After the  $\ell^{th}$-step,  the
digits $z_i$ will belong to the   alphabet
$$ \mathcal{A}^{(\ell)} = \{ -2a+\ell(B-M) ,  \ldots,   0, \ldots,  2a -\ell(B-M)\} \quad \hbox{for}\ 0 \leq \ell < s$$
and  to the alphabet $ \mathcal{A}$ for $ \ell = s$. Clearly
$\mathcal{A}^{(0)} = \{ -2a ,  \ldots, -1,  0,1, \ldots, 2a\}$.

\vskip0.3cm

\hrule

\vskip0.2cm

\noindent {\bf Algorithm~II}: Parallel addition for base $\beta$ with the weak representation of zero property ($\beta$ is WRZ).

\vskip0.2cm

\hrule

\vskip0.2cm

\noindent{\sl Input}:  two words $x_n \cdots x_m$ and $y_n \cdots
y_m$ of $\mathcal{A}^*$, with $m \le n$, $ x = \sum_{i=m}^{n}
x_i\beta^i$ and
$y = \sum_{i=m}^{n} y_i\beta^i$.\\
{\sl Output}: a word $z_{n+ks} \cdots z_{m-hs}$ of $\mathcal{A}^*$
such that
$$z =x+y = \sum_{i=m-hs}^{n+ks} z_i\beta^i.$$

\noindent\texttt{for each $i$ in parallel do}\\
0. \hspace*{0.5cm} $z_i := x_i+y_i$\\
1. \hspace*{0.5cm} \texttt{for} $\ell:=1$ \texttt{to} $s$ \texttt{do}\\
2. \hspace*{1.5cm} \texttt{if} $z_i \in \mathcal{A}'$
\texttt{then} $q_i :=0$ \texttt{else} $q_i := {\rm sgn}\, z_i$\\
3. \hspace*{1.5cm} $z_i:=z_i -\sum_{j=-h}^k  q_{i-j}b_j$

\vskip0.2cm

\hrule

\vskip0.5cm

\begin{thm}\label{AII}
Suppose that $\beta$ is WRZ.
Then Algorithm~II realizes addition in constant time in parallel in ${\rm Fin}_{\mathcal{A}}(\beta)$
with alphabet $\mathcal{A}= \{-a, \ldots,  0,\ldots,  a\}$, where $a=\left\lceil\tfrac{B-1}{2}\right\rceil + M$.
\end{thm}

\pf \noindent During the course of the algorithm, we do not change the value $\sum z_i\beta^i$
(thanks to the fact that we are only applying the weak representation of zero, either positively, or negatively).
For correctness of the algorithm, we just have to check the magnitude of digits at the time  when the algorithm stops.\\

After the $0^{th}$-step, any digit $z_i$ belongs to the alphabet $\mathcal{A}^{(0)}$.
We will  prove
\begin{description}
    \item[for ${1 \leq \ell < s}$\ :]  \quad  if after the $(\ell-1)^{th}$-step the digit $z_i$ belongs to $\mathcal{A}^{(\ell -1)}$, then
     after the $\ell^{th}$-step the new digit $z_i$  (we shall denote it for the moment by $z_i^{new}$) belongs to $ \mathcal{A}^{(\ell)}$.
    \item[for $\ell = s$\ :]  \quad  if after the $(s-1)^{th}$-step the digit $z_i$ belongs to  $\mathcal{A}^{(s -1)}$, then
     after the $s^{th}$-step the new digit $z_i^{new}$ belongs to $\mathcal{A}$.
\end{description}
This already will confirm the correctness of Algorithm~II.\\

Let us discuss  the value of digits $z_i$ after the $\ell^{th}$-step according to the value $q_i$  computed in this step.
During the discussion we shall  use  the inequalities
\begin{equation}\label{trebaII4} \begin{array}{lll}a< 2a-\ell
(B-M)&\hbox{for} & \ell <s\,,\\
a\geq 2a-\ell (B-M)&\hbox{for} & \ell =s\, ,
\end{array}
\end{equation}
which follow from the choice of $s$ by \eqref{trebaII3}. Let $z_i$ be in $\mathcal{A}^{(\ell-1)}$. There are three possible cases:

\begin{itemize}

    \item $q_i =0$: In this case, $|z_i|\leq a'=a -M$ and $z_i^{new} = z_i - \sum\limits_{j\neq 0}q_{i-j}b_j$.
     Therefore,
     $$ |z_i^{new}|  \leq |z_i| + |\sum_{j\neq 0}q_{i-j}b_j| \leq a-M+M  = a \leq \left\{\begin{array}{lll} 2a-\ell(B-M) &\hbox{if} & \ell < s\,,\\
     a &\hbox{if} & \ell = s\,\end{array}\right. $$
    {\it i.e.}, $z_i^{new} \in \mathcal{A}^{(\ell)}$ for $\ell < s$ and $z_i^{new} \in \mathcal{A}$ for  $\ell = s$.

    \item $q_i =1$: In this case, $ a'+1 = a-M+1\leq z_i \leq 2a -(\ell -1)(B-M) $ and $$z_i^{new} = z_i - B - \sum\limits_{j\neq 0}q_{i-j}b_j\,.$$
     Therefore, for the upper bound we have
     $$ z_i^{new}  \leq 2a-(\ell-1)(B-M)  -B+M = 2a-\ell(B-M)\leq \left\{\begin{array}{lll} 2a-\ell(B-M) &\hbox{if} & \ell < s\,,\\
    a &\hbox{if} & \ell = s\,. \end{array}\right. $$
    On the other hand, we obtain for the lower bound
    $$ z_i^{new} \geq a - M +1 -B - M \geq -a   \geq \left\{\begin{array}{lll} -2a+\ell(B-M) &\hbox{if} & \ell < s\,,\\
    -a &\hbox{if} & \ell = s\,, \end{array}\right.  $$
    {\it i.e.}, $z_i^{new} \in \mathcal{A}^{(\ell)}$ for $\ell < s$ and $z_i^{new} \in \mathcal{A}$ for  $\ell = s$.

    \item $q_i =- 1$: Analogous to the previous case.

\end{itemize}

In the previous discussion we have also  used the inequality $a - M +1 -B - M \geq -a$ which is a consequence  of the definition
of $a$, see \eqref{trebaII2}. \pfk

\begin{coro}
If $\beta$ is WRZ, then addition realized by Algorithm~II is a $(hs+ks+1)$-local function with memory $ks$ and anticipation $hs$. Algorithm~II is neighbour-free.
\end{coro}

\begin{remark}
If $\beta$ is a root of a strong polynomial, we may use both Algorithms~I and II with this strong representation of zero.
They work with the alphabet $\mathcal{A}= \{ -a, \ldots, 0, \ldots, a\}$, whereby
 in the first algorithm $a = a_I =\left\lceil\tfrac{B-1}{2}\right\rceil + \left\lceil\tfrac{B-1}{2(B-2M)}\right\rceil M\, $, while
 in the second one $a = a_{II} = \left\lceil\tfrac{B-1}{2}\right\rceil + M$.
Thus the alphabet of Algorithm~II is never bigger than the alphabet of Algorithm~I using the same SRZ.

It is easy to see that both the alphabets coincide ($a_I = a_{II}$) if, and only if,
$B\geq 4M-1$. If $B$ and $M$ satisfy this  inequality, then the parameter $s$ from
 Algorithm~II is equal to 1, and the  algorithms are the same.

If $4M-1 > B > 2M$, then Algorithm~II uses a strictly smaller alphabet.
The price we have to pay is the number of steps.
Unlike Algorithm~I, the number of steps in Algorithm~II depends on values $M$ and $B$ in the weak polynomial.
Another disadvantage of Algorithm~II is the increasing  number of positions needed to store the result.
Comparing the length of outputs in both algorithms, if $d=k+h$ is the degree of the strong or weak polynomial,
 Algorithm~I can enlarge the number of used positions by $d$, while Algorithm~II enlarges the number of used positions by $ds$.
\end{remark}

\begin{remark}
A given base $\beta$ can satisfy more than one strong polynomial, and, consequently, we have several versions of Algorithm~I. Similarly, one base $\beta$ can satisfy several different weak polynomials, and thus several versions of Algorithm~II.
\end{remark}

\begin{example} Let us  consider an integer base $\beta = b$ in $\mathbb{N}$.
\begin{itemize}
\item  If $b\geq 3$, then $-X + b = 0$ is a strong polynomial for $b$, with $B=b$ and $M=1$, thus $B\geq 4M-1$.
According to the previous remark, the two algorithms (I and II) coincide and require the same alphabet
$\mathcal{A} = \{ -a, \ldots, a \}$ with $a = \left\lceil\tfrac{b+1}{2}\right\rceil$.
In this sense, both algorithms are generalizations of the Avizienis algorithm to non-integer bases.

\item For $b=2$,  the equation   $-X + 2 = 0$ is a weak polynomial and we may use it in Algorithm~II,
 which works now with the alphabet $\mathcal{A} = \{-2, -1, 0, 1, 2\}$.
This alphabet is bigger than the alphabet in the parallel algorithm given by  Chow and Robertson, but,
 on the other hand, it is smaller than the alphabet in the algorithm that we have described in Example~\ref{binar}.
\end{itemize}
\end{example}

\begin{remark} In both our algorithms (I and II), the decision about application of the rewriting rule
at position $i$ depends only on the actual value of the digit at this position.
This is the crucial difference from the algorithm described by Chow and Robertson for base $2$, where the decision
depends also on the digit at the right neighboring position $i-1$.
A natural question is whether this kind of strategy may be applied to other bases as well. In the last section,
 we can answer this question positively, at least for the Golden Mean.
\end{remark}

%%%%%%%%%%%%%%%%%%%%%%%%%%%%%%%%%%%%%%%%%%%%%%%%%%%%%%%%%%%%%%%%%%%%%%%%%%%%%%%%%%%%%%%%%
\section{Polynomials with a dominant coefficient}

Algorithms I and II are applicable only to bases $\beta$ being a root of a strong or a weak polynomial, that is to say,
a polynomial with integer coefficients, where one of these coefficients is dominant,
\textit{i.e.}, is greater than twice or once the sum of the moduli of the other coefficients.
We shall show that for most of the algebraic numbers $\beta$ with $|\beta|>1$ such a strong or weak
polynomial exists.

\begin{prop}\label{dominant}
Let $\alpha$ be an algebraic number of degree $d$ with algebraic conjugates $\alpha_1, \alpha_2, \ldots, \alpha_d$
(including $\alpha$ itself). Let us assume that $|\alpha_i| \neq 1$ for all $i = 1,2,\ldots, d$ and  $|\alpha| >1$.
Then, for any $t \geq 1$ there exist a polynomial
$$ Q(X)= a_0 X^m + a_{1} X^{m-1} + \cdots + a_{m-1} X + a_m \in \mathbb{Z}[X] $$
and an index $i_0 \in \{ 1,\ldots, m\}$ such that
$$ Q(\alpha) = 0 \quad \hbox{and}\quad |a_{i_0}| > t \sum_{i\in \{0, \ldots, m\}\setminus \{i_0\}}{|a_i|}\,.$$
\end{prop}

\pf Let $G(X)$ be the minimal polynomial of $\alpha$, \textit{i.e.},
$$ G(X) = \prod_{j=1}^d (X-\alpha_j) = X^d + g_{1} X^{d-1} + \cdots + g_{d-1} X + g_{d} \in \mathbb{Q}[X]\,.$$
Let $M$ be the companion matrix of the polynomial $G(X)$:
$$M= \left(\begin{array}{llllll}
-g_{1} &1&0&0&\ldots &0\\
-g_{2} &0&1&0&\ldots &0\\
-g_{3} &0&0&1&\ldots &0\\
~~~~\vdots &&&&&\\
-g_{d-1} &0&0&0&\ldots &1\\
-g_{d} &0&0&0&\ldots &0
\end{array}\right) \in \mathbb{Q}^{d\times d}\,.$$
It is a well-known fact, which can be easily verified, that the characteristic polynomial of $M$ satisfies
$$ \det(M-XI)= (-1)^dG(X)\,, $$
where $I$ is the unit matrix of corresponding size $d \times d$. In particular, the numbers $\alpha_j$ are eigenvalues of $M$.
For any $n$ in $\mathbb{N}$,  $n\geq 2$, define
$$ G_n(X) = \prod_{j=1}^d (X-\alpha_j^n) = g_{0}(n) X^d + g_{1}(n) X^{d-1} + \ldots + g_{d-1}(n) X + g_{d}(n)\,.$$
Since  the matrix $M^n$ has eigenvalues $\alpha_1^n, \ldots, \alpha_d^n$, we have
$$ \det(M^n-XI) = (-1)^d\prod_{j=1}^d(X-\alpha_j^n) =  (-1)^d G_n(X)\,. $$
As $M$ is a rational  matrix, its power $M^n$  is rational as well, and thus $\det(M^n-XI)$ is a polynomial with rational coefficients.
It implies that, for all $n$ in $\mathbb{N}$, the polynomial $G_n(X)$ has rational coefficients.
Clearly, for all $n$ in $\mathbb{N}$ and $j$ in $\mathbb{N}, j \leq d$, we have
\begin{equation}\label{koef} g_j(n) = (-1)^j \sum_{\{i_1,i_2, \ldots, i_j\} \in S_j}
 \alpha_{i_1}^n\alpha_{i_2}^n\ldots \alpha_{i_j}^n\,,
\end{equation}
where $S_j = \bigl\{ A\subset \{1,2,\ldots,d\}\,: \, \# A = j \bigr\}$ is the set of all subsets of $\{ 1,2,\ldots, d\}$ with cardinality $j$.\\
Without loss of generality, let us assume $|\alpha_1| \geq  |\alpha_2| \geq \ldots \geq |\alpha_d|$, and
denote $j_0=\max\{ i : 1<|\alpha_i| \}$. Our choice of $j_0$ guarantees that
$$ \left|\frac{\alpha_{i_1} \alpha_{i_2} \ldots \alpha_{i_r}} { \alpha_1 \alpha_2 \ldots \alpha_{j_0} }\right| < 1 $$
for any subset $\{i_1,i_2, \ldots, i_r\} \subset \{1,2,\ldots, d\}$  and $\{i_1,i_2, \ldots, i_r\} \neq \{1,2,\ldots, j_0\}$.
Therefore, for all such choices of $\{ i_1, \ldots, i_r \} \neq \{1,2,\ldots, j_0\}$, we have
$$ \lim_{n\to \infty} \frac{\alpha_{i_1}^n \alpha_{i_2}^n \ldots \alpha_{i_r}^n} { \alpha_1^n \alpha_2^n \ldots \alpha_{j_0}^n }\ = \ 0  \,.$$
Since the coefficients $g_j(n)$ of the polynomial $G_n$ satisfy \eqref{koef},  we can deduce
$$ \lim_{n\to \infty}\,\frac{ g_{j}(n)}{~\alpha_1^n \alpha_2^n \ldots \alpha_{j_0}^n~} =
\left\{\begin{array}{cl} 0\quad & \hbox{for all}\ j = 1, \ldots, d, \, j\neq j_0 \,,\\
(-1)^j\quad & \hbox{for }\ j=j_0 \,.\\
\end{array}\right.$$
We may already claim that there exists  $n_0=n_0(t)$ in $\mathbb{N}$ such that
\begin{equation} \label{chceme1}
|g_{j_0}(n_0)| > t\sum_{j \in \{1, \ldots, m\}\setminus \{j_0\}} |g_{j}(n_0)|\,,
\end{equation}
or, equivalently,
$$ \frac{|g_{j_0}(n_0)|}{|\alpha_1^{n_0} \alpha_2^{n_0} \ldots \alpha_{j_0}^{n_0}|} >
 t\sum_{j \in \{1, \ldots, m\}\setminus \{j_0\}} \frac{|g_{j}(n_0)|}{|\alpha_1^{n_0} \alpha_2^{n_0} \ldots \alpha_{j_0}^{n_0}|}\,.$$
The existence of a suitable $n_0$ is obvious, as the right hand side tends to 0 and the left hand side has the limit $1$.

Let us fix one such $n_0$, and denote by $K$ the least common multiple of the denominators of ratios $g_{1}(n_0), \ldots, g_{d-1}(n_0), g_{d}(n_0)$.
Then the polynomial $Q(X)=KG_{n_0}(X^{n_0})$  and the index $i_0 = n_0 j_0$ have the desired properties.
\pfk

\begin{definition}\label{Def_t-pol}
Let $t \geq 1$, and let $T(X) = \sum_{i=0}^{n} t_{i} X^{i}$ be a polynomial with integer coefficients $t_{i}$ in $\mathbb{Z}$. $T$ is called a {\em $t$-polynomial} if there exists $i_0 \in \{ 0, \ldots, n \}$ such that
$$ |t_{i_0}| > t \sum_{i \in \{ 0, \ldots, n \} \setminus \{i_0\}} |t_{i}|\,.$$
\end{definition}

\begin{remark}
The polynomial $Q$ constructed in the proof of Proposition \ref{dominant} clearly is a~$t$-polynomial, and, for $t=1$ we obtain a weak polynomial, while for $t=2$ we get a strong polynomial.
\end{remark}

\begin{remark}  Lemma 8  in  \cite{Akiyama} gives a little bit weaker statement  in the case of  an expanding algebraic integer $\alpha$,
{\it i.e.}, $|\alpha_i| >1$ for all conjugates of $\alpha$. Our proof of the previous Proposition \ref{dominant} was strongly inspired  by the
proof given by S.~Akiyama, P.~Drungilas, and J.~Jankauskas in \cite{Akiyama}.
\end{remark}

\begin{example} The proof of Proposition \ref{dominant} is constructive. The strong representations of zero we have used in Examples
\ref{binar}, \ref{fibonacci}, and \ref{Gaussian} can be obtained by the construction given in the proof with $t=2$, namely as follows:
\begin{itemize}
    \item For $\beta = 2$, the minimal polynomial is $G(X) = X-2$, then $G_2(X) = X-2^2$, and the desired $Q(X) = G_2(X^2) = X^2 - 4$.
    \item For the Golden Mean, the minimal polynomial is $G(X) = X^2-X-1= (X-\beta)(X -\beta')$ with $\beta = \frac{1+\sqrt{5}}{2}$
        and $\beta' = \frac{1-\sqrt{5}}{2}$. The minimal $n$ satisfying \eqref{chceme1} is $n=4$. Therefore,
        $G_4(X) = \bigl(X-\beta^4\bigr)\bigl(X-(\beta')^4\bigr)$, and the final
        $Q(X) = G_4(X^4) = X^8 - \bigl(\beta^4 + (\beta')^4\bigr)X^4 + \beta^4 (\beta')^4 = X^8 -7X^4+1$.
    \item For the base $\beta = -1+i$, the minimal polynomial is $G(X) = (X+1-i)(X+1+i) = X^2 +2X +2$. Already for $n=2$ we have
        $G_2(X) = \bigl(X-(-1+i)^2\bigr)\bigl(X-(-1-i)^2\bigr) = (X+2i)(X-2i) = X^2+4$,
        and our strong representation of zero uses the polynomial $Q(X) = G_2(X^2) = X^4+4$.
\end{itemize}
\end{example}

\medskip

\begin{prop}\label{mod1}
Let $\beta$ be an algebraic number with a conjugate $\gamma$ such that $|\gamma|=1$.
Then $\beta$ cannot satisfy a $t$-polynomial for any $t \ge 1$.
\end{prop}
\pf
Suppose that $\beta$ is a root of a polynomial $W(X)=b_kX^k + b_{k-1}X^{k-1} + \cdots  + b_1X + b_0 + b_{-1}X^{-1} + \cdots + b_{-h}X^{-h}$ with integer coefficients $b_i \in \mathbb{Z}$. Then also the conjugate $\gamma$ is a root of $W(X)$:
$$b_k\gamma^k + b_{k-1}\gamma^{k-1} + \cdots  + b_1\gamma + b_0 + b_{-1}\gamma^{-1} + \cdots + b_{-h}\gamma^{-h}
=0.$$
Since $|\gamma|=1$, we obtain
\begin{equation}\label{eqmod1}
|b_0|\le \sum\limits_{i\in \{-h, \ldots,k\}\setminus \{0\}}{|b_i|}.
\end{equation}
\pfk

By Proposition~\ref{dominant} and Proposition~\ref{mod1} we obtain the following result:

\begin{thm}\label{SRR}
Let $\beta$  with $|\beta| >1$ be an algebraic number. Then $\beta$ is SRZ (or WRZ) if and only
if it has no conjugate of modulus $1$.
\end{thm}

\medskip

It is fairly easy to recognize whether an algebraic number does, or does not have a conjugate of modulus 1, by looking at its minimal polynomial.
First, if the number is quadratic, it cannot have any conjugate of modulus $1$.
Suppose now that $\alpha$ is an algebraic number of degree $d>2$, with a conjugate $\alpha'$ with modulus $|\alpha'| = 1$.  Let
$G(X) = X^d + g_{1}X^{d-1} + \cdots + g_{d-1}X +g_{d}$ be its minimal polynomial, $G(X)$ in $\mathbb{Q}[X]$.
Since $G(X)$ is minimal, $\alpha' \neq \pm 1$; {\it i.e.}, $\alpha'$ is not real. As the minimal polynomial has all its coefficients real,
the complex conjugate $\overline{\alpha'} = \frac{1}{\alpha'}$ is a root of $G$ as well.
In  general, if the minimal polynomial has two different roots $\eta$ and $\frac{1}{\eta}$, then the minimal polynomial satisfies $$ G(X) = X^d G\bigl(\tfrac{1}{X}\bigr)\, ,$$ thus it is reciprocal and its degree is even. This is summarized in the following remark.

\begin{remark}
Let $\beta$  with $|\beta| >1$ be an algebraic number of degree $d$.
\begin{itemize}
    \item If $d$ is odd, or
    \item if $d=2$, or
    \item if $d \geq 4$ is even and the minimal polynomial of $\beta$ is not reciprocal,
\end{itemize}
then $\beta$ has no conjugate of modulus $1$.
\end{remark}

%%%%%%%%%%%%%%%%%%%%%%%%%%%%%%%%%%%%%%%%%%%%%%%%%%%%%%%%%%%%%%%%%%%%%%%%%%%%%%%%%%%%%%%%%

\section{The Golden Mean - all good things come in threes}

In \cite{Berstel}, Berstel described an algorithm for parallel
addition in base $\beta= \frac{1+\sqrt{5}}{2}$ on an alphabet with
cardinality 13. In this section, we give Algorithm~III for
parallel addition in this numeration system on the alphabet $\{
-1,0, 1 \}$. This alphabet cannot be  further reduced, as proved
in \cite{Frougny2}. In our algorithm, we use a method similar to
the method of Chow and Robertson.

\bigskip

We begin by describing two auxiliary algorithms for elimination of digits. Both of them use the weak representation of zero
$-\beta^2 + 3 - \frac{1}{\beta^2} = 0$.

\bigskip

The first auxiliary algorithm removes digits $-2$:

\vskip0.2cm

\hrule

\vskip0.2cm

\noindent {\bf Algorithm~A}: Base $\beta= \frac{1+\sqrt{5}}{2}$, reduction from alphabet $\{ -2, \ldots, 2 \}$ to $\{ -1, \ldots, 2 \}$.

\vskip0.2cm

\hrule

\vskip0.2cm

\noindent{\sl Input}:  a finite sequence of digits $(z_i)$ of $\{-2,1,0,1,2\}$,
 with $ z= \sum z_i\beta^i$.\\
{\sl Output}: a finite sequence of digits $(z_i)$ of $\{-1
,0,1,2\}$,  with $ z= \sum z_i\beta^i$.

\vskip0.2cm

\noindent\texttt{for each $i$ in parallel do}\\
1. \hspace*{0.5cm} \texttt{case} $\left\{\begin{array}{l}
 z_i = -2\  \\
  z_i =  -1\  \\
 z_i =0\  \hbox{ \texttt{and}} \ z_{i+2} < 0\   \hbox{ \texttt{and}}\ z_{i-2}<   0\ \\
\end{array} \right\} $\ \texttt{then} $q_i:=-1$\\[2mm]
\hspace*{1cm} \texttt{else} $q_i:=0$\\
2. \hspace*{0.5cm} $z_i:=z_i-3q_i+q_{i+2}+q_{i-2}$

\vskip0.2cm

\hrule

\vskip0.5cm

\begin{proof}
For the correctness of Algorithm~A we have to show that the value $z_i^{new}=  z_i-3q_i+q_{i+2}+q_{i-2}$
belongs to the alphabet $\{-1,0,1,2\}$ for each $i$.\\
We will use the fact that our prescription for $q_i$ satisfies
\begin{equation}\label{removing2}
z_i\geq 0 \ \ {\rm and} \ \ z_{i+2}\geq 0 \qquad \Longrightarrow
\qquad q_i=0 \ \ {\rm and}\ \ q_{i+2}=0\,.
\end{equation}
Let us  discuss the value of $z_i^{new}$ according  to the values
$q_i$:
\begin{description}
\item[$q_i =0$] and
    \begin{itemize}
        \item $z_i\in\{1,2\}$: Then, $z_i^{new} = z_i + q_{i+2} + q_{i-2} \in \{-1,0,1,2\}$.
        \item $z_i=0$: Then our algorithm for the value $q_i$ supposes that $z_{i+2}\geq 0$ or $z_{i-2}\geq 0$. According to
         \eqref{removing2}, we have $q_{i+2}=0$ or $q_{i-2}=0$. Therefore, $z_i^{new} =q_{i+2} + q_{i-2} \in \{-1,0\}$.
    \end{itemize}
\item[$q_i=-1$] and
    \begin{itemize}
        \item $z_i\in\{-2,-1\}$: Then, $z_i^{new} = z_i + 3 + q_{i+2} + q_{i-2} \in \{-1, 0, 1, 2\}$.
        \item $z_i=0$: The rule for the value $q_i$ supposes that $z_{i+2} \leq -1$ and $z_{i-2} \leq -1$.
         Therefore, $q_{i+2} = q_{i-2} = -1$, which implies $z_i^{new} = 3 + q_{i+2} + q_{i-2} = 1$.
    \end{itemize}
\end{description}
Other combinations, namely ($q_i = 0$ and $z_i \in \{ -2, -1 \}$) or ($q_i = -1$ and $z_i \in \{ 1, 2 \}$), are impossible.
Thereby, we have shown that the digits of the result belong to the desired alphabet $\{ -1, 0, 1, 2\}$.
\end{proof}

\medskip

The second auxiliary algorithm removes digits $2$:

\vskip0.2cm

\hrule

\vskip0.2cm

\noindent {\bf Algorithm~B}: Base $\beta= \frac{1+\sqrt{5}}{2}$, reduction from alphabet $\{ -1, 0, 1, 2 \}$ to $\{ -1, 0, 1 \}$.

\vskip0.2cm

\hrule

\vskip0.2cm

\noindent{\sl Input}: a finite sequence of digits $(z_i)$ of $\{-1,0,1,2\}$,  with $ z= \sum z_i\beta^i$.\\
{\sl Output}: a finite sequence of digits $(z_i)$ of $\{-1, 0, 1\}$,  with $ z= \sum z_i\beta^i$.

\vskip0.2cm

\noindent\texttt{for each $i$ in parallel do}\\
1. \hspace*{0.5cm} \texttt{case} $\left\{\!\!\!\begin{array}{l}
 z_i = 2\  \\
 z_i =1\  \hbox{ \texttt{and}} \ (z_{i+2} \geq 1\   \hbox{ \texttt{or}}\ z_{i-2}\geq 1)\ \\
  z_i =0\  \hbox{ \texttt{and}} \ z_{i+2} =z_{i-2}=2\ \\
 z_i =0\  \hbox{ \texttt{and}} \ z_{i+2} = z_{i-2}=1\
 \hbox{ \texttt{and}} \ z_{i+4} \geq 1\   \hbox{ \texttt{and}}\ z_{i-4}\geq 1\\
 z_i =0\  \hbox{ \texttt{and}} \ z_{i+2} =2\   \hbox{ \texttt{and}}\ z_{i-2}=1\  \hbox{ \texttt{and}}\ z_{i-4}\geq 1\\
  z_i =0\  \hbox{ \texttt{and}} \ z_{i-2} =2\   \hbox{ \texttt{and}}\ z_{i+2}=1\  \hbox{ \texttt{and}}\ z_{i+4}\geq 1\\
\end{array}\!\!\! \right\} $\ \texttt{then} $q_i:=1$\\[2mm]
\hspace*{1cm} \texttt{else} $q_i:=0$\\
2. \hspace*{0.5cm} $z_i:=z_i-3q_i+q_{i+2}+q_{i-2}$

\vskip0.2cm

\hrule

\vskip0.5cm

\begin{proof}
Before explanation of the correctness of Algorithm~B, we point out several simple facts about the values $q_i$ determined in
Algorithm B:
\begin{itemize}
    \item Fact 1. \quad $z_i\geq 1 \ \ {\rm and} \ \ z_{i+2}\geq 1 \qquad \Longrightarrow \qquad q_i=1\ \ {\rm and}\ \ q_{i+2}=1\,.$
    \item Fact 2. \quad $z_i\leq 0 \ \ {\rm and} \ \ z_{i+2}\leq 0 \qquad \Longrightarrow \qquad q_i=0\ \ {\rm and}\ \ q_{i+2}=0\,.$
    \item Fact 3. \quad $z_{i+2}\leq 0 \ \ {\rm and} \ \ z_i=1 \ \ {\rm and} \ \ z_{i-2}\leq 0 \qquad \Longrightarrow \qquad q_i=0\,.$
    \item Fact 4. \quad $z_i= 0 \ \ {\rm and} \ \ q_i =1 \qquad \Longrightarrow \qquad q_{i+2}=1\ \ {\rm and}\ \ q_{i-2}=1\,.$
\end{itemize}

Unlike the previous four facts, the following three ones deserve short proofs:
\begin{itemize}
    \item Fact 5. \quad $z_i= 1 \ \ {\rm and} \ \ q_i =0 \qquad \Longrightarrow \qquad q_{i+2}=0\ \ {\rm and}\ \ q_{i-2}=0\,.$\\
        \textit{Proof of Fact 5.} The assumption $z_i= 1$ and  $q_i =0$ implies $z_{i+2}\leq 0$.
        If $z_{i+2} =-1$, then we have directly $q_{i+2}=0$. Now suppose that $z_{i+2}= 0$.
        If $q_{i+2}$ were equal to $1$, then, according to Fact 4, the value $q_i$ is equal
        to $1$, a contradiction.
        The discussion for $q_{i-2}$ is analogous.
    \item Fact 6. \quad $z_i= 0 \ \ {\rm and} \ \ q_i =0 \qquad \Longrightarrow \qquad q_{i+2}=0\ \ {\rm or}\ \ q_{i-2}=0\,.$\\
        \textit{Proof of Fact 6.} When $z_i= 0$,  then the  value $q_i$ is zero in  these two situations:\\
            a) One of numbers $z_{i+2}$, $z_{i-2}$ is not positive. Then, according to Fact 2, one of the values $q_{i+2}$, $q_{i-2}$ is zero.\\
            b) Both $z_{i+2}\geq 1 $ and $z_{i-2}\geq 1$.  The assumption $q_i=0$ forces, without loss of generality,
             that $z_{i-2}=1$ and  $z_{i-4}\leq 0$. So we can use Fact 3 to deduce $q_{i-2}=0$.
    \item Fact 7. \quad $z_i= 1 \ \ {\rm and} \ \ q_i =1 \qquad \Longrightarrow \qquad q_{i+2}=1\ \ {\rm or}\ \ q_{i-2}=1\,.$\\
        \textit{Proof of Fact 7.} The assumption $z_i= 1$ and $q_i =1$ implies $z_{i+2}\geq 1$ or $z_{i-2}\geq 1$.
        According to Fact 1, at least one of numbers $q_{i+2}$ or  $q_{i-2}$ is $1$.
\end{itemize}

Now we demonstrate the correctness of Algorithm~B, by showing that $z_i^{new}=  z_i-3q_i+q_{i+2}+q_{i-2}$
belongs to the alphabet $\{-1,0,1\}$ for each $i$.

\medskip

Let us  again discuss the value  $z_i^{new}$ according  to the values $q_i$:
\begin{description}
\item[$q_i =0$] and
\begin{itemize}
    \item $z_i=-1$: Then, $z_i^{new} =-1 + q_{i+2} + q_{i-2} \in \{-1,0,1\}$.
    \item $z_i=0$: Due to Fact 6, we have $z_i^{new} = q_{i+2}+q_{i-2} \in \{0,1\}$.
    \item $z_i=1$: According to Fact 5, we have  $z_i^{new} = 1+q_{i+2}+q_{i-2}=1$.
\end{itemize}
\item[$q_i=1$] and
\begin{itemize}
    \item $z_i=2$: Then, $z_i^{new} =-1 + q_{i+2} + q_{i-2} \in \{-1,0,1\}$.
    \item $z_i=1$: Due to Fact  7, we have  $z_i^{new} = -2+ q_{i+2}+q_{i-2} \in \{-1,0\}$.
    \item $z_i=0$: According to Fact 4, we have  $z_i^{new} = -3+q_{i+2}+q_{i-2}=-1$.
\end{itemize}
\end{description}
Other combinations, namely ($q_i = 0$ and $z_i =2$) or ($q_i = 1$ and $z_i = -1$), are impossible.
Thereby, we have shown that the digits of the result belong to the desired alphabet $\{ -1, 0, 1 \}$.
\end{proof}

\bigskip

Lastly, we can proceed by summarizing the Algorithm~III for parallel addition in base $\beta = \frac{1+\sqrt{5}}{2}$ and alphabet $\{ -1, 0, 1\}$:

\vskip0.2cm

\hrule

\vskip0.2cm

\noindent {\bf Algorithm~III}: Base $\beta= \frac{1+\sqrt{5}}{2}$, parallel addition on alphabet $\A=\{-1,0,1\}$.

\vskip0.2cm

\hrule

\vskip0.2cm

\noindent{\sl Input}:  two finite sequences of digits $(x_i)$ and $(y_i)$ of $\{-1,0,1\}$, with $ x = \sum x_i\beta^i$ and
$y = \sum y_i\beta^i$.\\
{\sl Output}: a finite sequence of digits $(z_i)$ of $\{-1,0,1\}$ such that $$z =x+y = \sum z_i\beta^i.$$

\noindent\texttt{for each $i$ in parallel do}\\
0. \hspace*{0.5cm} $v_i:=x_{i} +y_{i}$\\
1. \hspace*{0.5cm} use Algorithm~A
 with input
 $(v_i)$  and denote  its
 output  $(w_i)$\\
2. \hspace*{0.5cm} use Algorithm~B with input
  $(w_i)$ and denote  its
 output  $(z_i)$\\
\vskip0.2cm

\hrule

\vskip0.5cm

\begin{coro}\label{21}
Addition in base the Golden Mean on the alphabet $\{-1,0,1\}$ realized by Algorithm~III is a $21$-local function
with memory $10$ and anticipation $10$. Algorithm~III is neighbour-sensitive.
\end{coro}
\pf It is easily seen that Algorithm~A is local with memory $4$ and anticipation $4$.
Algorithm~B is local with memory $6$ and anticipation $6$.
Therefore, Algorithm~III is $21$-local with memory $10$ and anticipation $10$.
\pfk

\begin{remark}
All of our algorithms described above are working with symmetric alphabets of signed-digits. Such a choice of alphabet is practical, because the algorithm for addition can be used for subtraction at the same time. However, the symmetry of the alphabet is not necessary for parallelism as such. The parallel algorithm for
addition in base being the Golden Mean given by Berstel~\cite{Berstel} uses an alphabet of 13 non-negative digits, namely $\{0, 1 \ldots, 9, 10, 11, 12\}$. In fact, our auxiliary Algorithms A and B demonstrate that, if we do not insist on the symmetry, we may be able to further reduce the alphabet.
\end{remark}

\begin{remark}\label{2machines} In Algorithm~III, input values $x_{2i}$ and $y_{2i}$ at even positions do not influence the output values $z_{2i+1}$ at odd positions, and similarly for the input in odd positions. Therefore, the data at even and odd positions can be stored separately, and addition can be performed in parallel on the input $(x_{2i})$ and $(y_{2i})$ on one hand, and on the input $(x_{2i+1})$ and $(y_{2i+1})$ on the other hand, by an $11$-local function independently.
\end{remark}

\begin{remark}  The structure of Algorithm~III depends heavily on the coefficients of the weak polynomial satisfied by the Golden Mean. It is not clear how to generalize it to a broader class of bases. However, the case where $\beta$ is a quadratic unit integer seems feasible similarly to the Golden Mean.
\end{remark}

\begin{remark}{\bf (Fibonacci numeration system)} \ The Fibonacci numbers are defined as follows:
 $$F_0 = 1, \  F_1 = 2, \  \ \hbox{and} \ \ F_{n+2}= F_{n+1}+F_n \ \
 \hbox{for} \ n \in\mathbb{N}\,.$$
It is well known  that any non-negative  integer  $N$ can be expressed as a sum of different Fibonacci numbers,
$$N= F_{i_k} + F_{i_{k-1}} + \ldots + F_{i_0} \ \  \hbox{with}  \ \ i_k > i_{k-1} > \ldots > i_0\geq 0\,.$$
If, moreover, we require $i_j\geq i_{j-1}+2$ for any index $j$, then such an  expression is unique, see~\cite{Zeckendorf}.
For example, the number $29$ can be represented in the Fibonacci numeration system by $(29)_F = 1010000$, as $29 = F_6+F_4$.
If we abandon the requirement of uniqueness, and allow coefficients $-1,0,1$, then the numeration system is redundant, and addition can be performed by a parallel algorithm in constant time as well.
Let us explain it. We want to add  two integers  written in the form
$$X = \sum_{i\geq 0} x_{i}F_i\quad \hbox{and}  \quad Y = \sum_{i\geq 0} y_{i}F_i  \ \ \ \hbox{with}  \ \ x_i,y_i\in \{-1,0,1\}\,.$$
We can use Algorithm~III with a small modification.  Algorithm~III is based on the weak representation of zero  $-\beta^{i+2}+3 \beta^{i} -\beta^{i-2}
=0$. The Fibonacci numbers satisfy analogous relations:
\begin{itemize}
    \item $-F_{i+2} +3 F_{i} -F_{i-2} = 0$ for all $i\geq 2$
    \item $-F_{3} +3 F_{1} -F_{0} = 0$
    \item $-F_{2} +3 F_{0}  = 0$
\end{itemize}

\noindent Therefore, the Algorithm~III needs only slight changes, namely at the positions with indices $i=0, 1$ in Part 1 (Algorithm A) and
$i=0, 1, 2, 3$ in Part 2 (Algorithm B).

However, a function computable in parallel is defined by a sliding window of length $p$, which means that the function $\Phi$ cannot depend on any particular index $i$, but only on the $p$-tuples of the input alphabet.
To get rid of this problem, we use the classical trick of extending the alphabet $\{-1,0,1\}$ by an artificial symbol, say $ \$ $, to indicate the exceptional positions to the function $\Phi$.
Thus we will represent the above mentioned number $29$ as $(29)_F = 1010000\$\$\$\ldots$.
The domain of the $21$-local function $\Phi$ announced in Corollary \ref{21} has to be enlarged to the set $\{-1,0,1,\$\}^{21}$, and the definition of $\Phi$ must be correspondingly extended as well.

Note that, due to the relation $-F_{3} +3 F_{1} -F_{0} = 0$, the alternative option of processing the odd positions separately from the even positions (as described in Remark \ref{2machines}. for base $\beta = \frac{1+\sqrt{5}}{2}$) cannot be used for the Fibonacci numeration system.
\end{remark}

\section{Acknowledgements}
We are  grateful to Shigeki Akiyama for showing us the connection between polynomials with a dominant coefficient and the height reducing problem with expanding base $\alpha$, and for providing us the reference~\cite{Akiyama}. We also want to thank Wolfgang Steiner for pointing out that the results contained in~\cite{FrougnySteiner} imply Proposition~\ref{dominant} in case $\beta$ is a  Pisot number. Finally, we thank P\'eter Burcsi for Proposition~\ref{mod1}.

\medskip

We acknowledge financial support by the Czech Science Foundation grant GA\v{C}R 201/09/0584, and by the grants
MSM6840770039 and LC06002 of the Ministry of Education, Youth, and Sports of the Czech Republic.


\begin{thebibliography}{1}

\bibitem{Akiyama} \textsc{S. Akiyama, P. Drungilas, J. Jankauskas},  {Height reducing problem on algebraic integers}, preprint (2010), 22 pp.

\bibitem{Avizienis}  \textsc{A. Avizienis}, {Signed-digit number representations for fast parallel arithmetic}, {\em IRE Trans. Electron. Comput.} {\bf 10}
(1961) 389--400.

\bibitem{Berstel} \textsc{J. Berstel}, {Fibonacci words --- A survey}, {\em The  Book of L}, Springer-Verlag, 1986, 13--27.

\bibitem{Cauchy} \textsc{A. Cauchy}, Sur les moyens d'\'eviter les erreurs dans les calculs num\'eriques, {\em C.R. Acad. Sc. Paris} s\'erie I {\bf 11} (1840) 789--798.

\bibitem{ChowRobertson} \textsc{C.Y. Chow, J.E. Robertson}, {Logical design of a redundant binary adder}, in {\em Proc. 4th IEEE Symposium on Computer Arithmetic} (1978) 109--115.

\bibitem{Frougny2} \textsc{Ch. Frougny}, {On-line finite automata for addition in some numeration systems}, {\em RAIRO-Theor. Inf. Appl.} {\bf 33} (1999) 79--101.

\bibitem{FrougnyKlouda} \textsc{Ch. Frougny} and \textsc{K. Klouda}, Rational base number systems for $p$-adic numbers, {\em RAIRO-Theor. Inf. Appl.} (2011), to appear.

\bibitem{cant} \textsc{Ch. Frougny} and \textsc{J. Sakarovitch}, Number representation and finite automata, {\em Combinatorics, Automata and Number Theory}, V. Berth\'e, M. Rigo (Eds), Encyclopedia of Mathematics and its Applications 135, Cambridge University Press, 2010.

\bibitem{FrougnySteiner} \textsc{Ch. Frougny} and \textsc{W. Steiner},  {Minimal weight expansions in Pisot bases}, {\em Journal of Mathematical Cryptology} {\bf 2}  (2008) 365--392.

\bibitem{Kornerup} \textsc{P. Kornerup}, Necessary and Sufficient Conditions for Parallel, Constant Time Conversion and Addition, in {\em Proc. 14th IEEE Symposium on Computer Arithmetic} (1999) 152--155.

\bibitem{LM} \textsc{D. Lind} and \textsc{B. Marcus}, {\em An Introduction to Symbolic Dynamics and Coding}, Cambridge University Press, 1995.

\bibitem{Lothaire} \textsc{M. Lothaire}, \emph{Algebraic combinatorics on words}, Encyclopedia of Mathematics and its Applications 90, Cambridge University Press, 2002.

\bibitem{NM} \textsc{A.M. Nielsen} and \textsc{J.-M. Muller}, Borrow-Save Adders for Real and Complex Number Systems, in {\em Proc. Real Numbers and Computers}, Marseilles (1996) 121--137.

\bibitem{Penney} \textsc{W. Penney}, A ``binary" system for complex numbers. {\em Journal of the Association for Computing Machinery} {\bf 12} (1965) 247--248.

\bibitem{Zeckendorf} \textsc{E. Zeckendorf}, Repr\'esentation des nombres naturels par une somme de nombres de Fibonacci ou de nombres de Lucas, {\em  Bull. Soc. Roy. Sci. Li\`ege} {\bf 41} (1972) 179--182.

\end{thebibliography}
\end{document}